\documentclass[12pt,reqno]{amsart}
\usepackage{amssymb,a4wide,eucal,amscd,yhmath,amsgen,amsopn}
\usepackage{fancyhdr}
\usepackage{amsmath,hyperref,xcolor}
\usepackage{mathrsfs}
\usepackage[all,cmtip]{xy}
\SilentMatrices
\SelectTips{eu}{}
\newtheorem{theorem}{Theorem}[section]
\newtheorem{proposition}[theorem]{Proposition}

\newtheorem{lemma}[theorem]{Lemma}
\newtheorem{Not}[theorem]{Notation}
\theoremstyle{rem}

\theoremstyle{definition}
\newtheorem{definition}[theorem]{Definition}
\theoremstyle{construct}
\newtheorem{construct}[theorem]{Construction}
\theoremstyle{examp}

\newcommand\projective\mathbf
\newcommand\PP{\projective P}

\newcommand\OO{\mathcal O}

\newcommand\ZZ{\mathbb Z}

\newcommand\GG{\mathbb Gr}

\newcommand\onto\twoheadrightarrow
\newcommand\lra\longrightarrow
\newcommand\dar\downarrow
\DeclareMathOperator{\pic}{Pic}
\DeclareMathOperator{\im}{im}
\DeclareMathOperator{\cok}{coker}
\DeclareMathOperator{\rk}{rank}
\DeclareMathOperator{\Hom}{Hom}

\begin{document}

\title{Vector Bundle construction via Monads on multiprojective Spaces}
\author{Damian M Maingi}
\date{January, 2023}

\address{Department of Mathematics\\ School of Physical Sciences\\ Chiromo Campus\\ Chiromo Way\\University of Nairobi\\P.O Box 30197, 00100 Nairobi, Kenya\\
Department of Mathematics\\Sultan Qaboos University\\ P.O Box 50, 123 Muscat, Oman\\
Department of Mathematics\\Catholic University of Eastern Africa\\P.O Box 62157, 00200 Nairobi, Kenya}
\email{dmaingi@uonbi.ac.ke, dmaingi@cuea.edu, dmaingi@squ.edu.om}

\maketitle

\begin{abstract}
In this paper we construct indecomposable vector bundles associated to monads on multiprojective spaces.
Specifically we establish the existence of monads on $\PP^{2n+1}\times\PP^{2n+1}\times\cdots\times\PP^{2n+1}$ and on $\PP^{a_1}\times\cdots\times\PP^{a_n}$.
We prove stability of the kernel bundle which is a dual of a generalized Schwarzenberger bundle associated to the monads on $X=\PP^{2n+1}\times\PP^{2n+1}\times\cdots\times\PP^{2n+1}$
and prove that the cohomology vector bundle which is simple,  a generalization of special instanton bundles.
We also prove stability of the kernel bundle and that the cohomology vector bundle associated to the monad on $\PP^{a_1}\times\cdots\times\PP^{a_n}$ is simple.
Lastly, we construct explicitly the morphisms that establish the existence of monads  on $\PP^1\times\cdots\times\PP^1$.
\end{abstract}

\section{Introduction}
\noindent The existence of indecomposable low rank vector bundles on algebraic varieties in comparison with the ambient space
has been a fertile area in algebraic geometry for the last 45 years. In spite of this fact it remains intriguing and fascinating to construct new examples of indecomposable low rank vector bundles.
Some of the remarkable works in  this regard are: the famous Horrocks-Mumford bundle of rank 2 over $\PP^4$ \cite{11}, the Horrocks vector bundle of rank 3 on $\PP^5$ \cite{9}
the Tango bundles \cite{24} of rank $n-1$ on $\PP^n$ for $n\geq3$ and the rank 2 vector bundle on $\PP^5$ in characteristic 2 by Tango \cite{25} are all obtained as cohomologies of certain monads.\\
\\

\noindent This (monads) is one of the techniques used to construct these vector bundles. They appear in many contexts within algebraic geometry. 
and were first introduced by Horrocks\cite{10} where he proved that all vector bundles $E$ on $\PP^3$ could be obtained as the cohomology bundle of a given monad.
In vector bundle construction via monads on a given algebraic variety, the first task  is to show the existence of monads.
Fl\o{}ystad \cite{5} gave a theorem on the existence of monads over projective spaces. 
Costa and Miro-Roig \cite{3} extended these results to smooth quadric hypersurfaces of dimension at least 3.
Marchesi, Marques and Soares \cite{12} generalized Fl\o{}ystad's theorem to a larger set of varieties. 
Maingi \cite{16,17,18,19} proved the existence of monads on $\PP^n\times\PP^m$, $\PP^{2n+1}\times\PP^{2n+1}$,
$\PP^{a_1}\times\PP^{a_1}\times\PP^{a_2}\times\PP^{a_2}\times\cdots\times\PP^{a_n}\times\PP^{a_n}$ and on $\PP^n\times\PP^n\times\PP^m\times\PP^m$ respectively and proved simplicity of the cohomology bundles
associated.\\
\\
The flow of the results is similar to a paper by Ancona and Ottaviani \cite{1} where they proved that special instanton bundles on $\PP^{2n+1}$ are simple by 
first proving stability of Schwarzenberger bundles and also the results of Maingi \cite{18}.\\
\\
A natural and efficient technique to construct monads and hence more examples of vector bundles is to vary the ambient variety and
choose a different polarisation. In Section three of the paper, we first generalize the work of Maingi \cite{17} by construction of monads on $\PP^{2n+1}\times\cdots\times\PP^{2n+1}$ for a rank
$\beta-\alpha-\gamma$. We then prove stability of the kernel bundle which is a generalization of the dual of Schwarzenberger (steiner) bundles. Next we prove simplicity of the cohomolgy vector bundle.
Specifically we establish the existence of monads
\[\begin{CD}0@>>>{\OO_X(-1,\cdots,-1)^{\oplus\alpha}} @>>^{f}>\oplus\OO_X^{\oplus\beta}@>>^{g}>\OO_X(1,\cdots,1)^{\oplus\gamma} @>>>0\end{CD}\]
on $X=\PP^{2n+1}\times\cdots\times\PP^{2n+1}$. We shall call the monad above Type I in this paper.\\
\\
Next, in Section four we establish the existence of monads on $X=\PP^{a_1}\times\cdots\times\PP^{a_n}$ for the polarisation  
$\mathscr{L}=\OO_X(\alpha_1,\cdots,\alpha_t)$. This is a generalization of the results of Maingi \cite{16} Theorem 3.2 where he gave a conditional variant
theorem for the existence of a monad on $\PP^{n}\times\PP^{m}$.\\
\\
Specifically we establish the existence of monads
\[\begin{CD}0@>>>{\OO_X(-\alpha_1,\cdots,-\alpha_t)^{\oplus\alpha}} @>>^{f}>\oplus\OO_X^{\oplus\beta}@>>^{g}>\OO_X(\alpha_1,\cdots,\alpha_t)^{\oplus\gamma} @>>>0\end{CD}\]
on $X=\PP^{a_1}\times\cdots\times\PP^{a_n}$ which we shall call monad Type II.
We then prove stability of the kernel bundle $\ker g$ and finally prove that the cohomology vector bundle, $E=\ker g/\im f$ is simple.\\
\\
Lastly, in Section five we construct the morphisms that establish the existence of monads
\[\begin{CD}
M_\bullet: 0@>>>\OO_{X}(-1,\cdots,-1)^{\oplus{k}}@>>^{\overline{A}}>{\OO^{\oplus{2n}\oplus{2k}}_X} @>>^{\overline{B}}>\OO_{X}(1,\cdots,1)^{\oplus{k}}  @>>>0\\
\end{CD}\]
on $\PP^1\times\cdots\times\PP^1$ which are matrices whose entries are multidegree monomials.

\section{Preliminaries}

\noindent In this work we give generalizations for previous results by several authors. To be specific we build upon results by Maingi \cite{16,17,18} therefore the definitions, notation, the methods applied are quite similar and the trend follows the paper by 
Ancona and Ottaviani \cite{1}.In this section we define and give notation in order to set up for the main results.
Most of the definitions are from chapter two of the book by Okonek, Schneider and Spindler \cite{20}.

\begin{definition}
Let $X$ be a nonsingular projective variety. 
\begin{enumerate}
\renewcommand{\theenumi}{\alph{enumi}}
 \item A {\it{monad}} on $X$ is a complex of vector bundles:
\[
\xymatrix{0\ar[r] & M_0 \ar[r]^{\alpha} & M_1 \ar[r]^{\beta} & M_2 \ar[r] & 0}
\]
which is exact at $M_0$ and at $M_2$ i.e. $\alpha$ is injective and $\beta$ surjective.
\item A monad as defined above has a display diagram of short exact sequences as shown below:
\[
\begin{CD}
@.@.0@.0\\
@.@.@VVV@VVV\\
0@>>>{M_0} @>>>\ker{\beta}@>>>E@>>>0\\
@.||@.@VVV@VVV\\
0@>>>{M_0} @>>^{\alpha}>{M_1}@>>>\cok{\alpha}@>>>0\\
@.@.@V^{\beta}VV@VVV\\
@.@.{M_2}@={M_2}\\
@.@.@VVV@VVV\\
@.@.0@.0
\end{CD}
\]
\item The kernel of the map $\beta$, $F=\ker\beta$ and the cokernel of $\alpha$, $\cok\alpha$ for the given monad are also vector bundles and the vector bundle
$E = \ker(\beta)/\im (\alpha)$ and is called the cohomology bundle of the monad.
\end{enumerate}
\end{definition}

\begin{definition}
Let $X$ be a nonsingular projective variety, let $\mathscr{L}$ be a very ample line sheaf, and $V,W,U$ be finite dimensional $k$-vector spaces.
A linear monad on $X$ is a complex of sheaves,
\[ M_\bullet:
\xymatrix
{
0\ar[r] & V\otimes {\mathscr{L}}^{-1} \ar[r]^{A} & W\otimes \OO_X \ar[r]^{B} & U\otimes \mathscr{L} \ar[r] & 0
}
\]
where $A\in \Hom(V,W)\otimes H^0 \mathscr{L}$ is injective and $B\in \Hom(W,U)\otimes H^0 \mathscr{L}$ is surjective.\\
The existence of the monad $M_\bullet$ is equivalent to: $A$ and $B$ being of maximal rank and $BA$ being the zero matrix.
\end{definition}

\begin{definition}
Let $X$ be a non-singular irreducible projective variety of dimension $d$ and let $\mathscr{L}$ be an ample line bundle on $X$. For a 
torsion-free sheaf $F$ on $X$ we define
\begin{enumerate}
\renewcommand{\theenumi}{\alph{enumi}}
 \item the degree of $F$ relative to $\mathscr{L}$ as $\deg_{\mathscr{L}}F:= c_1(F)\cdot \mathscr{L}^{d-1}$, where $c_1(F)$ is the first Chern class of $F$
 \item the slope of $F$ as $\mu_{\mathscr{L}}(F):= \frac{\deg_{\mathscr{L}}F}{rk(F)}$.
  \end{enumerate}
\end{definition}

\subsection{Hoppe's Criterion over polycyclic varieties.}
Suppose that the Picard group Pic$(X) \simeq \ZZ^l$ where $l\geq2$ is an integer then $X$ is a polycyclic variety.
Given a divisor $B$ on $X$ we define $\delta_{\mathscr{L}}(B):= \deg_{\mathscr{L}}\OO_{X}(B)$.
Then one has the following stability criterion {\cite{13}, Theorem 3}:

\begin{theorem}[Generalized Hoppe Criterion]
 Let $G\rightarrow X$ be a holomorphic vector bundle of rank $r\geq2$ over a polycyclic variety $X$ equiped with a polarisation 
 $\mathscr{L}$ if
 \[H^0(X,(\wedge^sG)\otimes\OO_X(B))=0\] 
 for all $B\in\pic(X)$ and $s\in\{1,\ldots,r-1\}$ such that
 $\begin{CD}\displaystyle{\delta_{\mathscr{L}}(B)<-s\mu_{\mathscr{L}}(G)}\end{CD}$ then $G$ is stable and if
 $\begin{CD}\displaystyle{\delta_{\mathscr{L}}(B)\leq-s\mu_{\mathscr{L}}(G)}\end{CD}$ then $G$ is semi-stable.\\
\\
 Conversely if then $G$ is (semi-)stable then  \[H^0(X,G\otimes\OO_X(B))=0\]
 for all $B\in\pic(X)$ and all $s\in\{1,\ldots,r-1\}$ such that
 $\left(\delta_{\mathscr{L}}(B)\leq\right)$ $\delta_{\mathscr{L}}(B)<-s\mu_{\mathscr{L}}(G)$.
\end{theorem}

\vspace{0.75cm}

\begin{Not}
\noindent Suppose the ambient space is $X=\PP^{a_1}\times\cdots\times\PP^{a_n}$ then $\pic(X) \simeq \ZZ^{n}$.\\
We shall denote by $g_i$ for $i=1\cdots,n$ the generators of the Picard group of $X$, $\pic(X)$.\\
\\
Denote by $\OO_X(g_1, g_2,\cdots,g_{n}):= {p_1}^*\OO_{\PP^{a_1}}(g_1)\otimes\cdots\otimes {p_{n}}^*\OO_{\PP^{a_n}}(g_n)$,
where $p_i$ for $i=1,\cdots,n$ are natural projections  from $X$ onto $\PP^{a_i}$. \\
\\
For any line bundle $\mathscr{L} = \OO_X(g_1, g_2,\cdots,g_{n})$ on $X$ and a vector bundle $E$, we write 
$E(g_1, g_2,\cdots,g_{n}) = E\otimes\OO_X(g_1, g_2,\cdots,g_{n})$ 
and $(g_1, g_2,\cdots,g_{n}):= 1\cdot[g_1\times\PP^{a_1}]+\cdots+1\cdot[\PP^{a_n}\times g_{n}]$ representing its corresponding divisor.\\
\\
The normalization of $E$ on $X$ with respect to $\mathscr{L}$ is defined as follows:\\
Set $d=\deg_{\mathscr{L}}(\OO_X(1,0,\cdots,0))$, since $\deg_{\mathscr{L}}(E(-k_E,0,\cdots,0))=\deg_{\mathscr{L}}(E)-nk\cdot \rk(E)$ 
there is a unique integer $k_E:=\lceil\mu_\mathscr{L}(E)/d\rceil$ such that  $1 - d.\rk(E)\leq \deg_\mathscr{L}(E(-k_E,0,\cdots,0))\leq0$. 
The twisted bundle $E_{{\mathscr{L}}-norm}:= E(-k_E,0,\cdots,0)$ is called the $\mathscr{L}$-normalization of $E$.\\
\\
Lastly, the linear functional $\delta_{\mathscr{L}}$ on  $\mathbb{Z}^{n}$ is defined as $\delta_{\mathscr{L}}(p_1,p_2,\cdots,p_{n}):= \deg_{\mathscr{L}}\OO_{X}(p_1,p_2,\cdots,p_{n})$.\\
\\
For the $q^{-th}$ cohomology group we use the notation $H^q(\mathscr{F})$ in place of $H^q(X,\mathscr{F})$, for the sake of brevity.
\end{Not}

\noindent The following results are generalized versions to a multiprojective space for the purposes of this work.

\begin{proposition}
Let $X$ be a polycyclic variety with Picard number $n$, let $\mathscr{L}$ be an ample line bundle and
let E be a rank $r>1 $ holomorphic vector bundle over $X$.
If $H^0(X,(\bigwedge^q E)_{{\mathscr{L}}-norm}(p_1,\cdots,p_{n})) = 0$ for $1\leq q \leq r-1$ and every $(p_1,\cdots,p_{n})\in \mathbb{Z}^{n}$ such that $\delta_{\mathscr{L}}\leq0$
then E is $\mathscr{L}$-stable.
\end{proposition}

\begin{proposition}
Let $0\rightarrow E \rightarrow F \rightarrow G\rightarrow0$ be an exact sequence of vector bundles.
Then we have the following exact sequence involving exterior and symmetric powers
\[0\lra\bigwedge^q E \lra\bigwedge^q F \lra\bigwedge^{q-1} F\otimes G\lra\cdots \lra F\otimes S^{q-1}G \lra S^{q}G\lra0\]
\end{proposition}

\begin{theorem}[K\"{u}nneth formula]
 Let $X$ and $Y$ be projective varieties over a field $k$. 
 Let $\mathscr{F}$ and $\mathscr{G}$ be coherent sheaves on $X$ and $Y$ respectively.
 Let $\mathscr{F}\boxtimes\mathscr{G}$ denote $p_1^*(\mathscr{F})\otimes p_2^*(\mathscr{G})$\\
 then $\displaystyle{H^m(X\times Y,\mathscr{F}\boxtimes\mathscr{G}) \cong \bigoplus_{p+q=m} H^p(X,\mathscr{F})\otimes H^q(Y,\mathscr{G})}$.
\end{theorem}

\begin{lemma}
Let $X=\PP^{a_1}\times\cdots\times\PP^{a_n}$ then\\
$\displaystyle{H^t(X,\OO_X(p_1,\cdots,p_{n}) )\cong \bigoplus_{\sum_{q_i=1}^t} 
H^{q_1}(\PP^{a_1},\OO_{\PP^{a_1}}(p_1))\otimes H^{q_2}(\PP^{a_2},\OO_{\PP^{a_2}}(p_2))\otimes\cdots\otimes
H^{q_{n}}(\PP^{a_n},\OO_{\PP^{a_n}}(p_{n}))}$.
 \end{lemma}

\begin{theorem}[\cite{21}, Theorem 4.1]
 Let $n\geq1$ be an integer  and $d$ be an integer. We denote by $S_d$ the space of homogeneous polynomials of degree $d$ in 
 $n+1$ variables (conventionally if $d<0$ then $S_d=0$). Then the following statements are true:
 \begin{enumerate}
 \renewcommand{\theenumi}{\alph{enumi}}
  \item $H^0(\PP^n,\OO_{\PP^n}(d))=S_d$ for all $d$.
  \item $H^i(\PP^n,\OO_{\PP^n}(d))=0$ for $1<i<n$ and for all $d$.
  \item $H^n(\PP^n,\OO_{\PP^n}(d))\cong H^0(\PP^n,\OO_{\PP^n}(-d-n-1))$.
 \end{enumerate}
\end{theorem}

\begin{lemma}
If  $\displaystyle{\sum_{i=1}^np_i>}0$ then $h^p(X,\OO_X (-p_1,\cdots,-p_{n})^{\oplus k}) = 0$ where $X=\PP^{a_1}\times\cdots\times\PP^{a_n}$ and for $0\leq p< \dim(X) -1$, for $k$ a positive integer.
\end{lemma}

\begin{lemma}[\cite{12}, Lemma 10]
Let $A$ and $B$ be vector bundles canonically pulled back from $A'$ on $\PP^n$ and $B'$ on $\PP^m$ then\\
$\displaystyle{H^q(\bigwedge^s(A\otimes B))=
\sum_{k_1+\cdots+k_s=q}\big\{\bigoplus_{i=1}^{s}(\sum_{j=0}^s\sum_{m=0}^{k_i}H^m(\wedge^j(A))\otimes(H^{k_i-m}(\wedge^{s-j}(B)))) \big\}}$.
\end{lemma}

\begin{lemma}[\cite{5}, Main Theorem] Let $k\geq1$. There exists monads on $\PP^k$ whose maps are matrices of linear forms,
\[
\begin{CD}
0@>>>{\OO_{\PP^{k}}(-1)^{\oplus a}} @>>^{A}>{\OO^{\oplus b}_{\PP^{k}}} @>>^{B}>{\OO_{\PP^{k}}(1)^{\oplus c}} @>>>0\\
\end{CD}
\]
if and only if at least one of the following is fulfilled;\\
$(1)b\geq2c+k-1$ , $b\geq a+c$ and \\
$(2)b\geq a+c+k$
\end{lemma}

\begin{lemma}[\cite{17}, Theorem 3.9]
Let $n$ and $k$ be positive integers and $A$ and $B$ be morphisms of linear forms as in 
\[ B :=\left( \begin{array}{cccc|cccccccc}
x_0\cdots  & x_n &       &   &y_0 \cdots  & y_n\\
    &\ddots&\ddots &&\ddots&\ddots\\
    && x_0\cdots   x_n & & & y_0 \cdots  & y_n
\end{array} \right)
\]
 and
\[ A :=\left( \begin{array}{cccccccc}
-y_0\cdots  & -y_n \\
    	   &\ddots &\ddots\\
             &&-y_0 \cdots & -y_n\\
\hline
x_0 \cdots  & x_n \\
    	    &\ddots &\ddots\\
             && x_0\cdots & x_n\\
\end{array} \right)
\]
 then there exists a  linear monad of the form
\[
\begin{CD}
0@>>>\OO_{\PP^{2n+1}}(-1)^{\oplus k} @>>^{{A}}>{\OO^{\oplus2n+2k}_{\PP^{2n+1}}} @>>^{{B}}>\OO_{\PP^{2n+1}}(1)^{\oplus k} @>>>0\\
\end{CD}
\]
\end{lemma}

\begin{lemma}[\cite{16}, Theorem 3.2]
Let $X = \PP^n\times\PP^m$ and let $\mathscr{L} = \OO_X(\rho,\sigma)$ be an ample line bundle on $X$. Denote by $N = h^0(\OO_X(\rho,\sigma)) - 1$.
Let $\alpha,\beta, \gamma$ be positive integers such that at least one of the following conditions holds\\
 $(1)\beta\geq 2\gamma + N -1$, and $\beta\geq \alpha + \gamma$,  \\
$(2)\beta\geq \alpha + \gamma + N$. \\
\\
Then, there exists a linear monad on $X$ of the form
\[
\begin{CD}
0@>>>{\OO_{X}(-\rho,-\sigma)^{\oplus\alpha}} @>>^{A}>{\OO^{\oplus\beta}_X} @>>^{B}>{\OO_{X}(\rho,\sigma)^{\oplus\gamma}} @>>>0\\
\end{CD}
\]
\end{lemma}

\begin{definition}
Let $X$ be a projective variety. A sheaf $S$ on $X$ is a steiner bundle if has short exact sequence of the form
\[
\begin{CD}
0@>>>{\OO_{X}(-1)^{\oplus{a}}} @>>>{\OO^{\oplus{b}}_X} @>>>S @>>>0\\
\end{CD}
\]
They were first defined by Dolgachev and Kapranov \cite{4}.
\end{definition}

\begin{definition}\cite{23}
Let $k\geq0$ the exact sequence of sheaves on $\PP^{2n+1}$
\[
\begin{CD}
0@>>>{\OO_{X}(-1)^{\oplus{a}}} @>>^{\phi}>{\OO^{\oplus{b}}_X} @>>>S @>>>0\\
\end{CD}
\]
where $\phi$ is given by the matrix
\[ \left[ \begin{array}{cccc|cccccccc}
x_0\cdots  & x_n &       &   &y_0 \cdots  & y_n\\
    &\ddots&\ddots &&\ddots&\ddots\\
    && x_0\cdots   x_n & & & y_0 \cdots  & y_n
\end{array} \right]
\]
defines a $2n+k-$bundle $S$ on $\PP^{2n+1}$ called a (generalized) Schwarzenberger bundle. 
\end{definition}

\vspace{1cm}

\noindent The display of the monad in Lemma 2.14 is
\[
\begin{CD}
@.@.0@.0\\
@.@.@VVV@VVV\\
0@>>>{\OO_{\PP^{2n+1}}(-1)^{\oplus k}} @>>>S^*:=\ker{(B)}@>>>E@>>>0\\
@.||@.@VVV@VVV\\
0@>>>{\OO_{\PP^{2n+1}}(-1)^{\oplus k}} @>>^{A}>{\OO^{\oplus2n+2k}_{\PP^{2n+1}}}@>>>S:=\cok{(A)}@>>>0\\
@.@.@V^{B}VV@VVV\\
@.@.{\OO_{\PP^{2n+1}}(1)^{\oplus k}}@={\OO_{\PP^{2n+1}}(1)^{\oplus k}}\\
@.@.@VVV@VVV\\
@.@.0@.0
\end{CD}
\]
\noindent A special instanton bundle on $\PP^{2n+1}$ of quantum number $k$ is defined by the exact sequence
\[\begin{CD}
0@>>>{\OO_{\PP^{2n+1}}(-1)^{\oplus k}} @>>>S^*:=\ker{(B)}@>>>E@>>>0
\end{CD}\]
which is exactly the way Spindler and Trautmann remarkably described \cite{23} where $S$ is a Schwarzenberger bundle of rank $2n+k$ which 
is defined by the short exact sequence
\[\begin{CD}
0@>>>{\OO_{\PP^{2n+1}}(-1)^{\oplus k}} @>>^{A}>{\OO^{\oplus2n+2k}_{\PP^{2n+1}}}@>>>S:=\cok{(A)}@>>>0
\end{CD}\]
and they were proved by Ancona and Ottaviani \cite{1}, Theorem 2.2 to be stable and in Theorem 2.8 they proved that $E$ is simple. Independently, Bohnhorst and Spindler \cite{2} proved the stability
of rank $n$ Schwarzenberger bundles on $\PP^n$. 
In the next section we are going to establish the existence of monads on a more general space namely $\PP^{2n+1}\times\PP^{2n+1}\times\cdots\times\PP^{2n+1}$ and 
prove stability  of the kernel bundle $T$ and simplicity of the cohomology vector bundle $E$.

\vspace{1cm}

\section{Monad Type I and associated vector bundles}

\noindent The goal of this section is to construct monads over a multiprojective space of $m$ copies of $\PP^{2n+1}$.
More specifically we generalize the results of Maingi \cite{17} by varying the ambient space. We rely on methods similar to those used in \cite{18}.
The kernel bundle $T$ is a more generalized version of the dual of  a Schwarzenberger vector bundle and we prove that it is stable and consequently we prove that the cohomology
vector bundle $E$ associated to the monad on $X$ is simple. The vector bundle $E$ is a generalized version of an instanton bundle.

\begin{theorem}
Let $X=\PP^{2n+1}\times\cdots\times\PP^{2n+1}$ and $\mathscr{L} = \OO_X(1,\cdots,1)$ an ample line bundle. 
Denote by $N = h^0(\OO_X(1,\cdots,1)) - 1$. Then there exists a linear monad $M_\bullet$ on $X$ of the form
\[
\begin{CD}
M_\bullet: 0@>>>\OO_{X}(-1,\cdots,-1)^{\oplus{\alpha}}@>>^{f}>{\OO^{\oplus\beta}_X} @>>^{g}>\OO_{X}(1,\cdots,1)^{\oplus{\gamma}}  @>>>0\\
\end{CD}
\]
if atleast one of the following is satified
\begin{enumerate}
\renewcommand{\theenumi}{\alph{enumi}}
 \item $\beta\geq 2\gamma + N -1$, and $\beta\geq \alpha + \gamma$,
 \item $\beta\geq \alpha + \gamma + N$, where $\alpha,\beta, \gamma$ be positive integers. 
\end{enumerate} 
 
\end{theorem}

\begin{proof}
For the ample line bundle $\mathscr{L} = \OO_X(1,\ldots,1)$ we have the Segre embedding
\[
\xymatrix{
i^*:X = \PP^{2n+1}\times\cdots\times\PP^{2n+1} \hookrightarrow\PP\big(H^0(X,\OO_X(1,\ldots,1))\big)\cong \PP^{N:=(2n+2)^m-1}
}
\]
such that $i^*(\OO_X(1))\simeq \mathscr{L}$ 
\\
Suppose that one of the conditions of Lemma 2.13 is satified and we have $a=\alpha$, $b=\beta$, $c=\gamma$ and $k=2n+1$ thus there exists a linear monad
\[
\begin{CD}
0@>>>{\OO_{\PP^{2n+1}}(-1)^{\oplus\alpha}} @>>^{A}>{\OO^{\oplus\beta}_{\PP^{2n+1}}} @>>^{B}>{\OO_{\PP^{2n+1}}(1)^{\oplus\gamma}} @>>>0\\
\end{CD}
\]
on $\PP^{2n+1}$ whose morphisms are matrices $A$ and $B$ with entries monomials of degree one where
\begin{align*}
A\in\Hom(\OO_{\PP^{2n+1}}(-1)^{\oplus\alpha},\OO_{\PP^{2n+1}}^{\oplus\beta})\cong H^0(\PP^{2n+1},\OO_{\PP^{2n+1}}(1)^{\oplus\alpha\beta}) \\
B\in\Hom(\OO_{\PP^{2n+1}}^{\oplus\beta},\OO_{\PP^{2n+1}}(1)^{\oplus\gamma})\cong H^0(\PP^{2n+1},\OO_{\PP^{2n+1}}(1)^{\oplus\beta\gamma})
\end{align*}
Thus, $A$ and $B$ induce a monad on $X$, 
\[
\begin{CD}
0@>>>{{\mathscr{L}^{-1}}^{\oplus\alpha}} @>^{\bar{A}}>>{\OO^{\oplus\beta}_X} @>^{\bar{B}}>>{\mathscr{L}^{\oplus\gamma}} @>>>0\\
\end{CD}
\]
where 
whose morphisms are matrices $\bar{A}$ and $\bar{B}$ with entries multidegree monomials such that
\begin{center}
$\bar{A}\in\Hom(\OO_{X}(-\alpha_1,\ldots,-\alpha_t)^{\oplus\alpha},\OO^{\oplus\beta}_X)$ and
	$\bar{B}\in\Hom(\OO^{\oplus\beta}_X  , \OO_{X}(\alpha_1,\ldots,\alpha_t)^{\oplus\gamma})$
\end{center}
\end{proof}

\noindent The kernel bundle $T$ of the above monad is a generalization of the dual of Schwarzenberger vector bundles \cite{2} which we now 
proceed to prove that it is stable.

\begin{lemma}
Let $T$ be a vector bundle on $X=\PP^{2n+1}\times\cdots\times\PP^{2n+1}$ defined by the sequence
\[\begin{CD}0@>>>T @>>>\OO_X^{\oplus \beta}@>>>\OO_X(1,\cdots,1)^{\oplus\gamma} @>>>0\end{CD}\]
then $T$ is stable.
\end{lemma}

\begin{proof}
We show that $H^0(X,\bigwedge^q T(-p_1,\cdots,-p_m))=0$ for all $\displaystyle{\sum_i^m p_i>0}$ and $1\leq q\leq \rk(T)$.\\
\\
Consider the ample line bundle $\mathscr{L} = \OO_X(1,\cdots,1) = \OO(L)$. \\
Its class in 
$\pic(X)= \langle [g_1\times\PP^{2n+1}],\cdots,[\PP^{2n+1}\times g_m]\rangle$ corresponds to the class\\
$\displaystyle{\sum_{i=1}^m1\cdot[g_i\times\PP^{2n+1}]}$, where $g_i$, $i=1,\cdots,n$ are hyperplanes of $\PP^{2n+1}$ with the \\
intersection product induced by $g_i^{2n+1} = 1$ and $g_i^{2n+2}=0$.\\
\\
Now from the display diagram of the monad we get \\ 
\begin{align*}
\begin{split}
c_1(T) & = c_1(\OO_X^{\beta}) - c_1(\OO_X(1,\cdots,1)^{\oplus\gamma}) \\
       & = \beta(0,\cdots,0) - \gamma(1,\cdots,1)\\
       & = (-\gamma,\cdots,-\gamma)
\end{split}
\end{align*}
Now $L^{(2n+1)m}>0$ hence , the degree of $T$ is:
\begin{align*}
\begin{split}
\deg_{\mathscr{L}}T & = -\gamma([g_1\times\PP^{2n+1}]+\cdots+[\PP^{2n+1}\times g_m])\cdot (\sum_{i=1}^m1\cdot[g_i\times\PP^{2n+1}])^{m(2n+1)-1}\\
		    & = -\gamma L^{m(2n+1)}< 0.
\end{split}
\end{align*}
Since $\deg_{\mathscr{L}}T<0$, then $(\bigwedge^q T)_{\mathscr{L}-norm} = (\bigwedge^q T)$ and  it suffices by 
Proposition 2.6, to prove that $h^0(\bigwedge^q T(-p_1,\cdots,-p_m)) = 0$
with $\displaystyle{\sum_{i=1}^mp_i\geq0}$ and for all $1\leq q\leq \rk(T)-1$.\\
\\
Next we twist the exact sequence 
\[\begin{CD}
0@>>>T @>>>{\OO_X^{\oplus\beta}} @>>>\OO_X(1,\cdots,1)^{\oplus\gamma} @>>>0
\end{CD}\]
by $\OO_X(-p_1,\cdots,-p_m)$ we get,
\[
0\lra T(-p_1,\cdots,-p_m)\lra\OO_X(-p_1,\cdots,-p_m)^{\oplus\beta}\lra\OO_X(1-p_1,\cdots,1-p_m)^{\oplus\gamma}\lra0\]
and taking the exterior powers of the sequence by Proposition 2.7 we get
\[0\lra \bigwedge^q T(-p_1,\cdots,-p_m) \lra \bigwedge^q (\OO_X(-p_1,\cdots,-p_m)^{\oplus\beta})\lra \bigwedge^{q-1}(\OO_X(1-2p_1,\cdots,1-2p_m)^{\oplus\beta})\cdots\]
Taking cohomology we have the injection:
\[0\lra H^0(X,\bigwedge^{q}T(-p_1,\cdots,-p_m))\hookrightarrow H^0(X,\bigwedge^{q}(\OO_X(-p_1,\cdots,-p_m)^{\oplus\beta}))\]
Set $\mathscr{G}=\OO_X(-p_1,\cdots,-p_m)^{\beta} = \OO_X(-p_1,\cdots,-p_2)\otimes\OO_X^{\oplus\beta}$ and using Lemma 2.11
$H^0(X,\bigwedge^{q}\mathscr{G})$ expands into 
$\displaystyle{H^0(X,\sum_{j=0}^{q}\wedge^j\OO_X(-p_1,\cdots,-p_2)\otimes\OO_X^{\oplus\beta})}$ and since $\displaystyle{\sum_i^m p_i>0}$ by Lemma 2.12 then
\[h^0(X,\bigwedge^{q}(\OO_X(-p_1,\cdots,-p_m)^{\oplus\beta}))=h^0(X,\bigwedge^{q}T(-p_1,\cdots,-p_m))= 0\]
i.e. $h^0(\bigwedge^{q}T(-p_1,\cdots,-p_m))=0$ and thus $T$ is stable.

\end{proof}

\begin{theorem} Let $X={\PP^{2n+1}}\times\cdots\times{\PP^{2n+1}}$, then the cohomology vector bundle $E$ associated to the monad 
\[
\begin{CD}
0@>>>\OO_X(-1,\cdots,-1)^{\oplus\alpha} @>>^{{A}}>{\OO_X^{\oplus\beta}} @>>^{{B}}>\OO_X(1,\cdots,1)^{\oplus\gamma} @>>>0\\
\end{CD}
\]
of rank $\beta-\alpha-\gamma$ is simple.
\end{theorem}

\begin{proof}
The display of the monad is
\[
\begin{CD}
@.@.0@.0\\
@.@.@VVV@VVV\\
0@>>>{\OO_{X}(-1,\cdots,-1)^{\oplus\alpha}} @>>>T@>>>E@>>>0\\
@.||@.@VVV@VVV\\
0@>>>{\OO_{X}(-1,\cdots,-1)^{\oplus \alpha}} @>>^{f}>{\OO^{\oplus\beta}_{X}}@>>>Q@>>>0\\
@.@.@V^{g}VV@VVV\\
@.@.{\OO_{X}(1,\cdots,1)^{\oplus\gamma} }@={\OO_{X}(1,\cdots,1)^{\oplus\gamma} }\\
@.@.@VVV@VVV\\
@.@.0@.0
\end{CD}
\]

\noindent Since $E$ is simple if its only endomorphisms are the homotheties then we need to prove that $\Hom(E,E)=k$ which is equivalent to $h^0(E\otimes E^*)$.
\\
The first step is to take the dual short exact sequence \\
\[\begin{CD}
0@>>>\OO_X(-1,\cdots,-1)^{\oplus\alpha} @>>>T@>>>E @>>>0
\end{CD}\]
to get
\[
\begin{CD}
0@>>>E^* @>>>T^* @>>>\OO_X(1,\cdots,1)^{\oplus\alpha}@>>>0.
\end{CD}
\]
Tensoring by $E$ we get
\[
\begin{CD}
0@>>>E\otimes E^* @>>>E\otimes T^* @>>>E(1,\cdots,1)^{\oplus\alpha}@>>>0.
\end{CD}
\]
Now taking cohomology gives:
\[\begin{CD}
0@>>>H^0(X,E\otimes E^*) @>>>H^0(X,E\otimes T^*) @>>>H^0(E(1,\cdots,1)^{\oplus\alpha})@>>>\cdots
\end{CD}\]
which implies that 
\begin{equation}
h^0(X,E\otimes E^*) \leq h^0(X,E\otimes T^*)
\end{equation}
\\
Now we dualize the short exact sequence
\[\begin{CD}
0@>>>T @>>>{\OO_X^{\oplus\beta}} @>>>\OO_X(1,\cdots,1)^{\oplus\gamma} @>>>0
\end{CD}\]
to get
\[\begin{CD}
0@>>>\OO_X(-1,\cdots,-1)^{\oplus\gamma} @>>>{\OO_X^{\oplus\beta}} @>>>T^* @>>>0
\end{CD}
\]
Now twisting by $\OO_X(-1,\cdots,-1)$ and taking cohomology and get
\[\begin{CD}
0\lra H^0(X,\OO_X(-2,\cdots-2)^{\oplus\gamma}) \lra H^0(X,\OO_X(-1,\cdots,-1)^{\oplus\beta})\lra H^0(X,T^*(-1,\cdots,-1))\lra\\
\lra H^1(X,\OO_X(-2,\cdots,-2)^{\oplus\gamma}) \lra H^1(X,\OO_X(-1,\cdots,-1)^{\oplus\beta})\lra H^1(X,T^*(-1,\cdots,-1))\lra\\
\lra H^2(X,\OO_X(-2,\cdots,-2)^{\oplus\gamma}) \lra H^2(X,\OO_X(-1,\cdots,-1)^{\oplus\beta})\lra H^2(X,T^*(-1,\cdots,-1))\lra\cdots
\end{CD}
\]
from which we deduce $H^0(X,T^*(-1,\cdots,-1)) = 0$ and $H^1(X,T^*(-1,\cdots,-1)) = 0$ from Lemmas 2.9, 2.11 and Theorem 2.10.\\
\\
Lastly, tensor the short exact sequence
\[
\begin{CD}
0@>>>\OO(-1,\cdots,-1)^{\oplus\alpha} @>>>T @>>> E@>>>0\\
\end{CD}
\]
by $T^*$ to get
\[
\begin{CD}
0@>>>T^*(-1,\cdots,-1)^{\oplus\alpha} @>>>T\otimes T^* @>>> E\otimes T^*@>>>0\\
\end{CD}
\]
and taking cohomology we have
\[
\begin{CD}
0@>>>H^0(X,T^*(-1,\cdots,-1)^{\oplus\alpha}) @>>>H^0(X,T\otimes T^*) @>>> H^0(X,E\otimes T^*)@>>>\\
@>>>H^1(X,T^*(-1,\cdots,-1)^{\oplus\alpha})@>>>\cdots
\end{CD}
\]
But $H^1(X,T^*(-1,\cdots,-1)^{\oplus\alpha}=0$ for $\alpha>1$ from above.\\
\\
so we have 
\[
\begin{CD}
0@>>>H^0(X,T^*(-1,\cdots,-1)^{{\oplus\alpha}}) @>>>H^0(X,T\otimes T^*) @>>> H^0(X,E\otimes T^*)@>>>0
\end{CD}
\]
This implies that 
\begin{equation}
h^0(X,T\otimes T^*) \leq h^0(X,E\otimes T^*)
\end{equation}
\\
Since $T$ is stable then it follows that it is simple which implies $h^0(X,T\otimes T^*)=1$.\\
\\
From $(1)$ and now $(2)$ and putting these together we have;
\[1\leq h^0(X,E\otimes E^*) \leq h^0(X,E\otimes T^*) = h^0(X,T\otimes T^*) = 1\]
We have $ h^0(X,E\otimes E^*) = 1 $ and therefore $E$ is simple.

\end{proof}

\section{Monad Type II and associated vector bundles}

\noindent The goal of this section is to construct monads over a multiprojectivespace $\PP^{a_1}\times\cdots\times\PP^{a_n}$.
More specifically we generalize the results of Maingi \cite{16} by varying the ambient space and the polarisation $\mathscr{L}$.
We prove that the kernel bundle $F$ is stable and thereafter we prove that the cohomology
vector bundle $E$ associated to the monad on $X$ is simple.

\begin{theorem}
Let $X = \PP^{a_1}\cdots\times\PP^{a_n}$ and $\mathscr{L} = \OO_X(\alpha_1,\cdots,\alpha_t)$ an ample line bundle. 
Denote by $N = h^0(\OO_X(\alpha_1,\cdots,\alpha_t)) - 1$.
Then there exists a linear monad $M_\bullet$ on $X$ of the form
\[
\begin{CD}
M_\bullet: 0@>>>\OO_{X}(-\alpha_1,\cdots,-\alpha_t)^{\oplus\alpha}@>>^{f}>{\OO^{\oplus\beta}_X} @>>^{g}>\OO_{X}(\alpha_1,\cdots,\alpha_t)^{\oplus\gamma}  @>>>0\\
\end{CD}
\]
if atleast one of the following is satified
\begin{enumerate}
\renewcommand{\theenumi}{\alph{enumi}}
 \item $\beta\geq 2\gamma + N -1$, and $\beta\geq \alpha + \gamma$,
 \item $\beta\geq \alpha + \gamma + N$, where $\alpha,\beta, \gamma$ be positive integers. 
\end{enumerate}
\end{theorem}

\begin{proof}
For the ample line bundle $\mathscr{L} = \OO_X(\alpha_1,\ldots,\alpha_t)$ we have the Segre embedding
\[
\xymatrix{
i^*:X = \PP^{a_1}\cdots\times\PP^{a_n} \hookrightarrow\PP\big(H^0(X,\OO_X(\alpha_1,\ldots,\alpha_t))\big)\cong \PP^{N}
}
\]
such that $i^*(\OO_X(1))\simeq \mathscr{L}$ and
where $\displaystyle{N=\left({{a_1+\alpha_1}\choose\alpha_1}{{a_2+\alpha_2}\choose\alpha_2}\cdots{{a_n+\alpha_t}\choose\alpha_t}\right)-1}$
\\
Suppose that one of the conditions of Lemma 2.13 is satified thus there exists a linear monad
\[
\begin{CD}
0@>>>{\OO_{\PP^{N}}(-1)^{\oplus\alpha}} @>>^{A}>{\OO^{\oplus\beta}_{\PP^{N}}} @>>^{B}>{\OO_{\PP^{N}}(1)^{\oplus\gamma}} @>>>0\\
\end{CD}
\]
on $\PP^{N}$ whose morphisms are matrices $A$ and $B$ with entries monomials of degree one where
\begin{align*}
A\in\Hom(\OO_{\PP^{N}}(-1)^{\oplus\alpha},\OO_{\PP^{N}}^{\oplus\beta})\cong H^0(\PP^{N},\OO_{\PP^{N}}(1)^{\oplus\alpha\beta}) \\
B\in\Hom(\OO_{\PP^{N}}^{\oplus\beta},\OO_{\PP^{N}}(1)^{\oplus\gamma})\cong H^0(\PP^{N},\OO_{\PP^{N}}(1)^{\oplus\beta\gamma})
\end{align*}
Thus, $A$ and $B$ induce a monad on $X$, 
\[
\begin{CD}
0@>>>{{\mathscr{L}^{-1}}^{\oplus\alpha}} @>^{\bar{A}}>>{\OO^{\oplus\beta}_X} @>^{\bar{B}}>>{\mathscr{L}^{\oplus\gamma}} @>>>0\\
\end{CD}
\]
where 
whose morphisms are matrices $\bar{A}$ and $\bar{B}$ with entries multidegree monomials such that
\begin{center}
$\bar{A}\in\Hom(\OO_{X}(-\alpha_1,\ldots,-\alpha_t)^{\oplus\alpha},\OO^{\oplus\beta}_X)$ and
	$\bar{B}\in\Hom(\OO^{\oplus\beta}_X  , \OO_{X}(\alpha_1,\ldots,\alpha_t)^{\oplus\gamma})$
\end{center}
\end{proof}

\begin{theorem}
Let $F$ be a vector bundle on $X =\PP^{a_1}\times\cdots\times\PP^{a_n}$ defined by the short exact sequence
\[\begin{CD}0@>>>F @>>>\OO_X^{\oplus\beta}@>>^{g}>\OO_X(\alpha_1,\cdots,\alpha_t)^{\oplus\gamma} @>>>0\end{CD}\]
then $F$ is stable for an ample line bundle $\mathscr{L} = \OO_X(\alpha_1,\cdots,\alpha_t)$
\end{theorem}

\begin{proof}
We are going to show that $H^0(X,\bigwedge^q F(-p_1,\cdots,-p_{n}))=0$ for all $\displaystyle{\sum_{i=1}^{n}p_i>0}$ and $1\leq q\leq \rk(F)-1$.\\
\\
Consider the ample line bundle $\mathscr{L} = \OO_X(\alpha_1,\cdots,\alpha_t) = \OO(L)$. \\
Its class in $\pic(X)= \langle [g_i\times\PP^{a_i}],i=1,\ldots,n]\rangle$ corresponds to 
$\displaystyle{\sum_{i=1}^{n}1.[g_i\times\PP^{a_i}]}$
where each $g_i$ is a hyperplane in $\PP^{a_i}$ with intersection product induced by 
$g_i^{a_i} =1$ and $g_i^{a_i+1}=0$ for $i=1,\ldots,n$.\\ 
\\
From the display of the monad we get
\[c_1(F) = c_1(\OO_X^{\oplus\beta}) - c_1(\OO_X(\alpha_1,\cdots,\alpha_t)^{\oplus\gamma}) =(-\gamma\alpha_1,\cdots,-\gamma\alpha_t)\]
Since $L^{a_1+\cdots+a_n}>0$ the degree of $F$ is $\deg_{\mathscr{L}}F = c_1(T)\cdot\mathscr{L}^{d-1}$ that is\\
\begin{align*}
\begin{split}
=-\gamma{n}{\sum_{i=1}^t\alpha_i}([g_1\times\PP^{a_1}]+\cdots+[\PP^{a_n}\times g_{2n}])\left(\displaystyle{\sum_{i=1}^n1\cdot[g_i\times\PP^{a_i}]}\right)^{{\sum_{i=1}^{n}a_i}-1}\\
=-\gamma{n}{\sum_{i=1}^t\alpha_i}L^{(a_1+\cdots+a_n)}< 0
\end{split}
\end{align*}
\\
Since $\deg_{\mathscr{L}}F<0$, then $(\bigwedge^q F)_{\mathscr{L}-norm} = (\bigwedge^q F)$ and  it suffices by 
the generalized Hoppe Criterion (Proposition 2.6), to prove that $h^0(\bigwedge^q F(-p_1,-p_2,\cdots,-p_{n})) = 0$
with $\displaystyle{\sum_{i=1}^{n}p_i>0}$ and for all $1\leq q\leq \rk(F)-1$.\\
\\
Next consider the exact sequence 
\[\begin{CD}0@>>>F @>>>\OO_X^{\oplus\beta}@>>^{g}>\OO_X(\alpha_1,\cdots,\alpha_t)^{\oplus\gamma} @>>>0\end{CD}\]
on twisting it by $\OO_X(-p_1,\cdots,-p_{n})$ one gets,
\[\begin{CD}0@>>>F(-p_1,\cdots,-p_{n}) @>>>\OO_X^{\oplus\beta}(-p_1,\cdots,-p_{n})@>>^{g}>\OO_X(\alpha_1-p_1,\cdots,\alpha_t-p_n)^{\oplus\gamma} @>>>0\end{CD}\]
and taking the exterior powers of the sequence by Proposition 2.10 one gets
\[0\lra \bigwedge^q F(-p_1,\cdots,-p_{n}) \lra \bigwedge^q (\OO_X(-p_1,\cdots,-p_{n})^{\oplus\beta})\lra \bigwedge^{q-1} (\OO_X(\alpha_1-2p_1,\cdots,\alpha_t-2p_{n})^{\oplus\gamma})\lra \cdots\]
Taking cohomology we have the injection:
\[0\lra H^0(X,\bigwedge^q F(-p_1,\cdots,-p_{n}))\hookrightarrow H^0(X,\bigwedge^q (\OO_X(-p_1,\cdots,-p_{n})^{\oplus\beta})\]
From here $h^0(X,\bigwedge^{q}F(-p_1,\cdots,-p_n)) = 0$ is proved in the same way as Lemma 3.2 the last part and thus $F$ is stable.

\end{proof}

\begin{theorem} Let $X =\PP^{a_1}\times\cdots\times\PP^{a_n}$, then the cohomology vector bundle $E$ associated to the monad 
\[\begin{CD}0@>>>{\OO_X(-\alpha_1,\cdots,-\alpha_t)^{\oplus\alpha}} @>>^{f}>\OO_X^{\oplus\beta}@>>^{g}>\OO_X(\alpha_1,\cdots,\alpha_t)^{\oplus\gamma} @>>>0\end{CD}\]
of rank $\beta-\alpha-\gamma$ is simple.
\end{theorem}

\begin{proof}
The display of the monad is
\[
\begin{CD}
@.@.0@.0\\
@.@.@VVV@VVV\\
0@>>>{\OO_{X}(-\alpha_1,\cdots,-\alpha_t)^{\oplus\alpha}} @>>>F=\ker g@>>>E@>>>0\\
@.||@.@VVV@VVV\\
0@>>>{\OO_{X}(-\alpha_1,\cdots,-\alpha_t)^{\oplus\alpha}} @>>^{f}>{\OO_X^{\oplus\beta}}@>>>Q=\cok f@>>>0\\
@.@.@V^{g}VV@VVV\\
@.@.{\OO_{X}(\alpha_1,\cdots,\alpha_t)^{\oplus\gamma}}@={\OO_{X}(\alpha_1,\cdots,\alpha_t)^{\oplus\gamma} }\\
@.@.@VVV@VVV\\
@.@.0@.0
\end{CD}
\]

\noindent Since $T$ is stable from Lemma 4.2 we prove that the cohomology vector bundle $E$ with rank $2n$ is simple.\\
\noindent On taking the dual of the short exact sequence on the first row of the display diagram and tensoring by $E$ we obtain
\[
\begin{CD}
0@>>>E\otimes E^* @>>>E\otimes F^* @>>>E(t,\cdots,t)^{\oplus\alpha}@>>>0.
\end{CD}
\]
Now taking cohomology gives:
\[\begin{CD}
0@>>>H^0(X,E\otimes E^*) @>>>H^0(X,E\otimes F^*) @>>>H^0(E(\alpha_1,\cdots,\alpha_t)^{\oplus\alpha})@>>>\cdots
\end{CD}\]
\\
which implies that 
\begin{equation}
h^0(X,E\otimes E^*) \leq h^0(X,E\otimes F^*)
\end{equation}
\\
Dualize the short exact sequence on the first column of the display diagram to get
\[\begin{CD}
0@>>>\OO_X(-\alpha_1,\cdots,-\alpha_t)^{\oplus\gamma} @>>>{\OO_X^{\beta}} @>>>F^* @>>>0
\end{CD}\]
\\
Now twisting the short exact sequence above by $\OO_X(-\alpha_1,\cdots,-\alpha_t)$ one obtains the short exact sequence
\[\begin{CD}
0@>>>\OO_X(-2\alpha_1,\cdots,-2\alpha_t)^{\oplus\gamma} @>>>{\OO_X(-\alpha_1,\cdots,-\alpha_t)^{\beta}} @>>>F^*(-\alpha_1,\cdots,-\alpha_t) @>>>0
\end{CD}\]
next on taking cohomology one gets\\
\[\begin{CD}
0\lra H^0(\OO_X(-2\alpha_1,\cdots,-2\alpha_t)^{\oplus\gamma}) \lra H^0(\OO_X(-\alpha_1,\cdots,-\alpha_t)^{\beta})\lra H^0(F^*(-\alpha_1,\cdots,-\alpha_t))\lra\\
0\lra H^1(\OO_X(-2\alpha_1,\cdots,-2\alpha_t)^{\oplus\gamma}) \lra H^1(\OO_X(-\alpha_1,\cdots,-\alpha_t)^{\beta})\lra H^1(F^*(-\alpha_1,\cdots,-\alpha_t))\lra\\
\lra H^2(\OO_X(-2\alpha_1,\cdots,-2\alpha_t)^{\oplus\gamma}) \lra H^2(\OO_X(-\alpha_1,\cdots,-\alpha_t)^{\beta})\lra H^2(F^*(-\alpha_1,\cdots,-\alpha_t))\lra\cdots
\end{CD}
\]
\\
from which we deduce $H^0(X,F^*(-\alpha_1,\cdots,-\alpha_t)) = 0$ and $H^1(X,F^*(-\alpha_1,\cdots,-\alpha_t)) = 0$ from Lemmas 2.9, 2.11 and Theorem 2.10.\\
\\
Lastly, tensor the short exact sequence
\[
\begin{CD}
0@>>>\OO(-\alpha_1,\cdots,-\alpha_t)^{\oplus k} @>>>F @>>> E@>>>0\\
\end{CD}
\]
by $F^*$ to get
\[
\begin{CD}
0@>>>F^*(-\alpha_1,\cdots,-\alpha_t)^k @>>>F\otimes F^* @>>> E\otimes F^*@>>>0\\
\end{CD}
\]
and taking cohomology we have
\[
\begin{CD}
0@>>>H^0(X,F^*(-\alpha_1,\cdots,-\alpha_t)^k) @>>>H^0(X,F\otimes F^*) @>>> H^0(X,E\otimes F^*)@>>>\\
@>>>H^1(X,F^*(-\alpha_1,\cdots,-\alpha_t)^k)@>>>\cdots
\end{CD}
\]
\\
But since  $H^0(X,F^*(-\alpha_1,\cdots,-\alpha_t)) = H^1(X,F^*(-\alpha_1,\cdots,-\alpha_t)) = 0$ from above then  it follows $H^1(X,F^*(-\alpha_1,\cdots,-\alpha_t)^k)=0$ for $k>1$.\\
\\
so we have 
\[
\begin{CD}
0@>>>H^0(X,F^*(-\alpha_1,\cdots,-\alpha_t)^{k}) @>>>H^0(X,F\otimes F^*) @>>> H^0(X,E\otimes F^*)@>>>0
\end{CD}
\]
This implies that 
\begin{equation}
h^0(X,F\otimes F^*) \leq h^0(X,E\otimes F^*)
\end{equation}
\\
Since $F$ is stable then it is simple implying $h^0(X,F\otimes F^*)=1$.\\
\\
From $(3)$ and $(4)$ and putting these together we have;\\
\[1\leq h^0(X,E\otimes E^*) \leq h^0(X,E\otimes F^*) = h^0(X,F\otimes F^*) = 1\]\\
\\
We have $ h^0(X,E\otimes E^*) = 1 $ and therefore $E$ is simple.

\end{proof}

\section{Monad construction via morphisms}
\noindent Let $X$ be a nonsingular projective variety. A monad$\xymatrix{0\ar[r] & M_0 \ar[r]^{\alpha} & M_1 \ar[r]^{\beta} & M_2 \ar[r] & 0}$ on $X$ exists
if one can give the morphisms $\alpha$ and $\beta$. In this section we establish the existence of monads on $\PP^1\times\cdots\times\PP^1$ by providing an explicit contruction of the morphisms
derived from the matrices used by Fl\o{}ystad \cite{5} and Ancona and Ottaviani \cite{1}.

\vspace{1cm}

\begin{construct}
Let $\psi : X = \PP^1\times\cdots\times\PP^1\longrightarrow \PP^{N=2n+1}$ be the Segre embedding which is defined as follows:\\
\\
$[\alpha_{10}:\alpha_{11}][\alpha_{20}:\alpha_{21}]:\ldots:[\alpha_{m0}:\alpha_{m1}]\hookrightarrow [x_0:x_1:\cdots:x_n:y_0:y_2:\ldots:y_n]$.\\
\\
First note that since we are taking $m$ copies of $\PP^1$ then we have \\
$N=2^m-1$\\
$=2^m-2+1$\\
$=2(2^{m-1}-1)+1$\\
$=2n+1$\\
\\i.e. $N=2n+1$ where $m$ and $n$ are positive integers such that $n=2^{m-1}-1$.\\
\\
Thus from Lemma 2.14, there exists a linear monad
\[\begin{CD}
0@>>>{\OO_{\PP^{2n+1}}(-1)^{\oplus{k}}} @>>^{A}>{\OO^{\oplus2n+2k}_{\PP^{2n+1}}}@>>^{B}>\OO_{\PP^{2n+1}}(1)^{\oplus{k}}@>>>0\end{CD}
\]
whose morphisms $A$ and $B$ that establish the monad are as given in Lemma 2.14. \\
\\
We induce a monad on $X=\PP^1\times\cdots\times\PP^1$%
 \[\begin{CD}
M_\bullet: 0@>>>\OO_{X}(-1,\cdots,-1)^{\oplus{k}}@>>^{\overline{A}}>{\OO^{\oplus{2n}\oplus{2k}}_X} @>>^{\overline{B}}>\OO_{X}(1,\cdots,1)^{\oplus{k}}  @>>>0\\
\end{CD}\]
by giving the morphisms $\overline{A}$ and $\overline{B}$ with $\overline{B}\cdot\overline{A}=0$ and $\overline{A}$ and $\overline{B}$ are of maximal rank.\\
\\%
From $A$ and $B$ whose entries are $x_0,\cdots,x_n,y_0,\cdots,y_n$ the homogeneous coordinates on $\PP^{2n+1}$ we give the correspondence
for the the Segre embedding using the following table:\\
\[ \begin{array}{|c|c|}
\hline
homog. coord. ~~on ~~\PP^{2n+1} & representation ~homog. coord.~~ on ~~X\\
\hline
x_0 & a_{0000\cdots0000} \\
x_1 & a_{0000\cdots0001} \\
x_2 & a_{0000\cdots0010} \\
x_3 & a_{0000\cdots0011} \\
x_4 & a_{0000\cdots0100} \\
\vdots&\vdots\\
x_{n-1} &  a_{0111\cdots1110} \\
x_n &  a_{0111\cdots1111} \\
y_0 & a_{1000\cdots0000} \\
y_1 & a_{1000\cdots0001} \\
y_2 & a_{1000\cdots0010} \\
y_3 & a_{1000\cdots0011} \\
y_4 & a_{1000\cdots0100} \\
\vdots&\vdots\\
y_{n-1} & a_{1111\cdots1110} \\
y_n & a_{1111\cdots1111} \\
\hline
\end{array} 
\]
\\
where $a_{iiii\cdots iiii}$ for $i$ is $0$ or $1$ are monomials of multidegree $(1,\ldots,1)$ i.e.\\

\[ \begin{array}{|c|c|}
\hline
representation ~~homog. coord. ~~on ~~\PP^{2n+1} & homog. coord.~~ on ~~X\\
\hline
a_{0000\cdots0000} & \alpha_{10}\alpha_{20}\alpha_{30}\alpha_{40}\cdots\alpha_{(m-3)0}\alpha_{(m-2)0}\alpha_{(m-1)0}\alpha_{m0} \\
a_{0000\cdots0001} & \alpha_{10}\alpha_{20}\alpha_{30}\alpha_{40}\cdots\alpha_{(m-3)0}\alpha_{(m-2)0}\alpha_{(m-1)0}\alpha_{m1} \\
a_{0000\cdots0010} & \alpha_{10}\alpha_{20}\alpha_{30}\alpha_{40}\cdots\alpha_{(m-3)0}\alpha_{(m-2)0}\alpha_{(m-1)1}\alpha_{m0} \\
a_{0000\cdots0011} & \alpha_{10}\alpha_{20}\alpha_{30}\alpha_{40}\cdots\alpha_{(m-3)0}\alpha_{(m-2)0}\alpha_{(m-1)1}\alpha_{m1} \\
a_{0000\cdots0100} & \alpha_{10}\alpha_{20}\alpha_{30}\alpha_{40}\cdots\alpha_{(m-3)0}\alpha_{(m-2)1}\alpha_{(m-1)0}\alpha_{m0} \\
\vdots&\vdots\\
a_{0111\cdots1110} & \alpha_{10}\alpha_{21}\alpha_{31}\alpha_{41}\cdots\alpha_{(m-3)1}\alpha_{(m-2)1}\alpha_{(m-1)1}\alpha_{m0} \\
a_{0111\cdots1111} & \alpha_{10}\alpha_{21}\alpha_{31}\alpha_{41}\cdots\alpha_{(m-3)1}\alpha_{(m-2)1}\alpha_{(m-1)1}\alpha_{m1} \\
a_{1000\cdots0000} & \alpha_{11}\alpha_{20}\alpha_{30}\alpha_{40}\cdots\alpha_{(m-3)0}\alpha_{(m-2)0}\alpha_{(m-1)0}\alpha_{m0} \\
a_{1000\cdots0001} & \alpha_{11}\alpha_{20}\alpha_{30}\alpha_{40}\cdots\alpha_{(m-3)0}\alpha_{(m-2)0}\alpha_{(m-1)0}\alpha_{m1} \\
a_{1000\cdots0010} & \alpha_{11}\alpha_{20}\alpha_{30}\alpha_{40}\cdots\alpha_{(m-3)0}\alpha_{(m-2)0}\alpha_{(m-1)1}\alpha_{m0} \\
a_{1000\cdots0011} & \alpha_{11}\alpha_{20}\alpha_{30}\alpha_{40}\cdots\alpha_{(m-3)0}\alpha_{(m-2)0}\alpha_{(m-1)1}\alpha_{m1} \\
a_{1000\cdots0100} & \alpha_{11}\alpha_{20}\alpha_{30}\alpha_{40}\cdots\alpha_{(m-3)0}\alpha_{(m-2)1}\alpha_{(m-1)0}\alpha_{m0} \\
\vdots&\vdots\\
a_{1111\cdots1110} & \alpha_{11}\alpha_{21}\alpha_{31}\alpha_{41}\cdots\alpha_{(m-3)1}\alpha_{(m-2)1}\alpha_{(m-1)1}\alpha_{m0} \\
a_{1111\cdots1111} & \alpha_{11}\alpha_{21}\alpha_{31}\alpha_{41}\cdots\alpha_{(m-3)1}\alpha_{(m-2)1}\alpha_{(m-1)1}\alpha_{m1} \\
\hline
\end{array} 
\]
Specifically we define $\overline{A}$ and $\overline{B}$ as follows\\
\[ \overline{B} :=\left[ \begin{array}{ccc|cccccc}
a_{0000\cdots0000}~~~\cdots  & a_{0111\cdots1111} &  &a_{1000\cdots0000} ~~~\cdots  & a_{1111\cdots1111}\\
    \ddots&\ddots &&\ddots&\ddots\\
    & a_{0000\cdots0000}~~~ \cdots~~~  a_{0111\cdots1111}  && & a_{1000\cdots0000}~~~ \cdots  & a_{1111\cdots1111}
\end{array} \right]
\]
\\
 and
\[ \overline{A} :=\left[ \begin{array}{cccccccc}
-a_{1000\cdots0000}\cdots  & -a_{1111\cdots1111} \\
    	   &\ddots &\ddots\\
             &&-a_{1000\cdots0000} \cdots & -a_{1111\cdots1111}\\
\hline
a_{0000\cdots0000}\cdots  & a_{0111\cdots1111} \\
    	    &\ddots &\ddots\\
             && a_{0000\cdots0000}\cdots & a_{0111\cdots1111}\\
\end{array} \right]
\]

We note that
\begin{enumerate}
 \item $\overline{B}\cdot \overline{A} = 0$ and
 \\
 \item  The matrices $\overline{B}$ and $\overline{A}$ have maximal rank\\
\end{enumerate}
Hence we get the desired monad,
\[\begin{CD}
M_\bullet: 0@>>>\OO_{X}(-1,\cdots,-1)^{\oplus{k}}@>>^{\overline{A}}>{\OO^{\oplus{2n}\oplus{2k}}_X} @>>^{\overline{B}}>\OO_{X}(1,\cdots,1)^{\oplus{k}}  @>>>0\\
\end{CD}\]
\end{construct}

\section{Acknowledgment}
\noindent I wish to express sincere thanks to the Department of Mathematics and Actuarial Science, Catholic University of Eastern Africa and the Department of Mathematics, University of Nairobi, for providing 
a conducive enviroment to be able to carry out research despite the overwhelming duties in teaching and community service.
I am also extremely grateful to Melissa, my wife who always encourages me to keep on and lastly to Amelia, Jerome and Wachuka who are always around when I am working on my research
at home.

\vspace{1cm}

\noindent \textbf{Data Availability statement}
My manuscript has no associate data.

\vspace{1cm}

\noindent \textbf{Conflict of interest}
On behalf of all authors, the corresponding author states that there is no conflict of interest.


\vspace{2cm}

\newpage

\begin{thebibliography}{1}
\bibitem[1]{1} Ancona V and Ottaviani G: \textit{Stability of special instanton Bundles on $\PP^{2n+1}$}.
Transactions of the American Mathematical Society 341 (1994) 677 - 693. 
\href{https://doi.org/10.2307/2154578}{\color{blue}doi: 10.2307/2154578.}
%
\bibitem[2]{2} Bohnhorst G and Spindler H (1992). The stability of certain vector bundles on $\PP^{n}$. 
Lecture Notes in Mathematics. Vol 1507, Springer, Berlin, Heidelberg. \href{https://doi.org/10.1007/BFb0094509}{\color{blue}doi: 10.1007/BFb0094509}
%
\bibitem[3]{3} Costa L and Mir\'{o}-Roig R M. Monads and instanton bundles on smooth hyperquadrics. Mathematische Nachrichten, 282 (2009), no 2, 169-179. 
\href{https://onlinelibrary.wiley.com/doi/10.1002/mana.200610730}{\color{blue}doi: 10.1002/mana.2000610730.}
%

\bibitem[4]{4} Dolgachev I and Kapranov M. Arrangements of hyperplanes and vector bundles on $\PP^n$.  Duke Math. J. 71 (1993), 633-664.
 \href{https://doi.org/10.1215/S0012-7094-93-07125-6}{\color{blue}doi: 10.1215/S0012-7094-93-07125-6}
%
\bibitem[5]{5} Fl\o{}ystad G. Monads on a Projective Space. Communications in Algebra, 28 (2000), 5503 - 5516.
\href{https://www.tandfonline.com/doi/abs/10.1080/00927870008827171}{\color{blue}doi: 10.1080/00927870008827171.}
%
\bibitem[6]{6} Hartshorne R. Varieties of small codimension in projective space.  Bull. Amer. Math. Soc. 80 (1974), 1017-1032.
\href{https://doi.org/10.1090/S0002-9904-1974-13612-8}{doi 10.1090/S0002-9904-1974-13612-8}
%
\bibitem[7]{7} Hartshorne R. Algebraic vector bundles on projective spaces: A problem list, Topology, 1979, vol 18, 117-128
\href{https://www.sciencedirect.com/science/article/pii/0040938379900302?via%3Dihub}{\color{blue}doi.org/10.1016/0040-9383(79)90030-2}
%

\bibitem[8]{8} Hoppe H. Generischer Spaltungstyp und zweite Chernklasse stabiler Vektorraumbündel. vom Rang 4 auf $\PP^4$ , Math. Z. 187 (1984), 345–360.
\href{http://gdz.sub.uni-goettingen.de/dms/resolveppn/?PPN=GDZPPN002429659}{\color{blue}eudml.org/doc/173496.}
%

\bibitem[9]{9} Horrocks G. Examples of rank three vector bundles on five-dimensional projective space,
Journal of London Mathematical Society. (2), 18 (1978), 15-27.
 \href{https://londmathsoc.onlinelibrary.wiley.com/doi/epdf/10.1112/jlms/s2-18.1.15}{\color{blue}doi:10.1112/jlms/s2-18.1.15.}

%

\bibitem[10]{10} Horrocks G. Vector bundles on the punctured spectrum of a local ring. Proc. London Math.
Soc. 14, 689-713 (1964). \href{https://londmathsoc.onlinelibrary.wiley.com/doi/abs/10.1112/plms/s3-14.4.689}{\color{blue}doi.org/10.1112/plms/s3-14.4.689.}

%
\bibitem[11]{11} Horrocks G and Mumford D. A rank 2 vector bundle on $\PP^{4}$ with 15,000 symmetries, Topology, 12 (1973), 63-81. 
\href{https://doi.org/10.1016/0040-9383(73)90022-0}{\color{blue}doi:10.1016/0040-9383(73) 90022-0.}
%
\bibitem[12]{12} Jardim M and Earp HNS\'{a}. Monad construction of asymptotically stable Bundles. arXiv (September,2011).
\href{https://arxiv.org/pdf/1109.2750v1.pdf}{\color{blue}arxiv.org/pdf/1109.2750v1.pdf}

%
\bibitem[13]{13} Jardim M, Menet M, Prata D and Earp HNS\'{a}. Holomorphic bundles for higher dimensional gauge theory.
Bulletin London Mathematical society, 49 (2017). \href{https://doi.org/10.1112/blms.12017}{\color{blue}doi: 10.1112/blms.12017.}
%


%
\bibitem[14]{14} Kumar, N. Construction of rank two vector bundles on $\PP^4$ in positive characteristic. 
Invent math 130, 277–286 (1997). \href{https://doi.org/10.1007/s002220050185}{\color{blue}doi 10.1007/s002220050185}


\bibitem[15]{15} Kumar N, Peterson C and Rao A P. Construction of low rank vector bundles on $\PP^4$ and $\PP^5$, Journal of Algebraic Geometry 11 (2) (2002), 203–217.
\href{https://doi.org/10.1090/S1056-3911-01-00309-5}{\color{blue}doi 10.1090/S1056-3911-01-00309-5}


\bibitem[16]{16} Maingi D. Vector Bundles of low rank on a multiprojective space. Le Matematiche. Vol. LXIX (2014) - Fasc. II. pp 31-41.
\href{https://lematematiche.dmi.unict.it/index.php/lematematiche/article/view/1048}{\color{blue}doi: 10.4418/2014.69.2.4.}
%
\bibitem[17]{17} Maingi D (2021). Indecomposable Vector Bundles associated to Monads on Cartesian products of projective spaces.
Turkish Journal of Mathematics. Vol. 45: No. 5. Article 17. Pages 2126-2139.
\href{https://doi.org/10.3906/mat-2101-6}{\color{blue}doi: 10.3906/mat-2101-6}

%
\bibitem[18]{18} Maingi D. Monads on multiprojective Products of Projective Spaces.
Manuscripta Mathematica (2022),\href{https://doi.org/10.1007/s00229-022-01449-0}{\color{blue}doi: 10.1007/s00229-022-01449-0}.

\bibitem[19]{19} Maingi D. Vector Bundles associated to monads on Cartesian Products of Projective Spaces.
Open Journal of Mathematical Sciences, OMS - Vol 7 (2023), Issue 1, pp 148-159\href{https://doi.org/10.30538/oms2023.0203}{\color{blue}doi: 10.30538/oms2023-0203}.


\bibitem[20]{20} Marchesi S, Marques P M and Soares H. Monads on a Projective Varieties. Pacific Journal of Mathematics, vol 296 (2018), no. 1, 155-180. 
\href{https://doi.org/10.2140/pjm.2018.296.155}{\color{blue}doi: 10.2140/pjm.2018.296.155.}
%

\bibitem[21]{21} Okonek C, Schneider M and Spindler H. Vector Bundles on Complex Projective Spaces.
Springer, 1980, \href{https://doi.org/10.1007/978-1-4757-1460-9}{\color{blue}doi.org/10.1007/978-1-4757-1460-9}
%
\bibitem[22]{22} Perrin D. G\'{e}om\'{e}trie alg\'{e}brique. Une introduction, (1995), EDP Sciences/CNRS \'{e}dition.
\href{https://www.decitre.fr/livres/geometrie-algebrique-9782729605636.html}{\color{blue}ISBN-2-7296-0563-0}
%

\bibitem[23]{23} Spindler H and Trautmann G. Special Instanton bundles on $\PP^{2n+1}$ their geometry and their moduli. Mathematische Annalen 286. 1-3 (1990): 559-592.
 \href{https://doi.org/10.1215/S0012-7094-93-07125-6}{\color{blue}doi: 10.1215/S0012-7094-93-07125-6}


\bibitem[24]{24} Tango, H. An example of indecomposable vector bundle of rank $n-1$ on $\PP^n$, $n\geq3$. 
Journal of Mathematics of Kyoto University, 16, (1976): 137-141,
\href{http://doi.org/10.1215/kjm/1250522965}{\color{blue}doi: 10.1215/kjm/1250522965.}

%
\bibitem[25]{25} Tango H. On morphisms from projective space $\PP^{n}$ to the Grassmann variety $\GG(n,d)$,
Journal of Mathematics of Kyoto University, 16 (1976), 201-207. 
\href{http://doi.org/10.1215/kjm/1250522969}{\color{blue}doi: 10.1215/kjm/1250522969.}

\end{thebibliography}
\end{document}